\documentclass[reqno]{amsart}
\usepackage{amsfonts}
\usepackage{amsmath}
\usepackage{graphics}
\usepackage{mathrsfs}
\usepackage{ dsfont}
\newcommand{\nr}{{|\!|}}
\newcommand{\st}{\int_0^t}
\newcommand{\s}{\int_{\Omega}}
\newcommand{\R}{{\mathbb{R}^3}}

\newtheorem{thm}{Theorem}[section]

\newtheorem{remark}[thm]{Remark}

\newcommand{\rem}[1]{}

\newcommand{\RR}{{\mathbb{R}^3}}
\newcommand {\sys}{$\mathcal S$}

\def\dd{\delta}

\begin{document}

\title[Anisotropic 3D Brinkman-Forchheimer-extended Darcy Model]
{Existence and uniqueness of global solutions for the modified anisotropic 3D Navier-Stokes equations
}
\date{}

\author{Hakima Bessaih}
\address[H. Bessaih] {University of Wyoming, Department of Mathematics, Dept. 3036, 1000 East University
Avenue, Laramie WY 82071, United States.}
\email{bessaih@uwyo.edu}
\author[S. Trabelsi]{Saber Trabelsi}
\address[S. Trabelsi]
{ Division of Math and Computer Sci. and Eng. \\
 King Abdullah University of Science and Technology \\
Thuwal 23955-6900 \\
Saudi Arabia
}
\email{saber.trabelsi@kaust.edu.sa}
\author{Hamdi Zorgati}
\address[H. Zorgati] {D\'epartement de Math\'ematiques, Campus Universitaire, Universit\'e Tunis El Manar 2092, Tunisia.}
\email{hamdizorgati@yahoo.fr}

\begin{abstract}
We study a modified three-dimensional incompressible anisotropic Navier-Stokes equations. The modification consists in the addition of a power term to the nonlinear convective one. This modification appears naturally in porous media when a fluid obeys the Darcy-Forchheimer law instead of the classical Darcy law. We prove global in time existence and uniqueness of solutions  without assuming the smallness condition on the initial data. This improves the result obtained for the classical $3D$ incompressible anisotropic Navier-Stokes equations. 
\end{abstract}
\maketitle

\section*{Introduction}
The purpose of this paper is to study the following modified Navier-Stokes system
\begin{align*}
\mathcal S_a:\quad \left\lbrace\begin{array}{rcll}
\partial_t\,u-\nu\,\Delta_h\,u +(u\cdot\nabla)\,u +a\,|u|^{2\alpha}u&=&-\nabla\,p\quad &\text{ for}\quad (t,x)\in\R\times \RR,\\ \\
 \nabla\cdot u&=&0 \quad &{\text{ for}}\quad (t,x)\in \R\times \RR,\\ \\
u_{\vert_{t=0}}&=&u_0,& 
\end{array}
\right.
\end{align*}
\vskip6pt
\noindent 
where $\partial_t$ denotes the partial derivative with respect to time, $\alpha\in \mathbb R, a>0, \Delta_h:=\partial_1^2 +\partial_2^2$ and $\partial_i$ denotes the partial derivative  in the direction $x_i$. Clearly, 
$\mathcal S_0$ corresponds to the classical anisotropic Navier-Stokes equations.   When a Coriolis force 
$\frac1\epsilon u\times {\bf e}_3$ is added, where ${\bf e}_3$ denotes the unit vertical vector and $\epsilon>0$ is the so called Rossby number, system $\mathcal S_0$ models rotating flows (see e.g. \cite{Pedo}). 
 We refer the reader to e.g. \cite{Nader} and \cite{porous} where the relevance of considering 
anisotropic viscosities of the form $\nu_h\Delta_h u+\epsilon\nu_v\,\partial_3^2u$ is explained through 
the Ekman's law. For a complete discussion leading to the anisotropic Navier-Stokes systems, the reader is referred to the book \cite{Pedo} or the introduction of the book \cite{Chemin0}. System $\mathcal S_a$ with $\Delta_h$ replaced by the classical Laplacian  is nothing but the three-dimensional Brinkman-Forchheimer-extended Darcy model. 
The Brinkman-Forchheimer-extended Darcy  equations 
\begin{equation}\label{ddbf}\partial_t\,u-\nu\Delta u+(u\cdot\nabla)\,u +\nabla\,p+a\,|u|^{2\alpha}u=f, \quad  
\nabla\cdot u=0,\end{equation}
have been extensively studied. The existence of weak solutions for $\alpha\geq 0$ and existence (for $\alpha\geq\frac54$) and uniqueness (for $\frac54<\alpha\leq 2$) of strong solutions of \eqref{ddbf} is shown in \cite{Cai}. Also, in \cite{MTT}, existence and uniqueness of weak and  strong solutions is shown for a larger range of $\alpha$. In \cite{Varga}, the authors show existence and uniqueness of solutions for all $\alpha>1$, with Dirichlet boundary conditions and regular enough initial data. Their argument relies on the maximal regularity estimate for the corresponding semi-linear stationary Stokes problem proved using some modification of the nonlinear localization technique. In \cite{MTT}, the authors showed the existence and uniqueness of weak and strong solutions, in particular with initial data in $H^1$ instead of $H^2$ in \cite{Varga} with periodic boundary conditions.  Let us mention that the space $L^{2\alpha+2}$ appears naturally in the mathematical analysis of system \eqref{ddbf}, and it coincides obviously with $L^4$ for $\alpha=1$. Since $\dot H^\frac12$ is known to be the critical Sobolev space for the classical 3D Navier-Stokes equations and $\dot H^\frac12\subset L^4$, then it makes sense to assume $\alpha>1$. For a more detailed discussion about  equation \eqref{ddbf} and the various values of $\alpha$ that lead to  the well posedness, we refer the reader to  \cite{Varga,MTT} and the references therein.

\vskip6pt
\noindent
Before going further, let us precise the notation and the functional setting that will be used along the paper. Since the horizontal variable $x_h\dd(x_1,x_2)$ does not play the same role as the vertical variable $x_3$, it is natural to introduce functional spaces taking into account this feature. These spaces are the so called anisotropic Sobolev spaces $H^{s, s'}$ for all $s,s' \in \mathbb R$.  More precisely, the space $H^{s, s'}$  is the Sobolev space with regularity $H^s$ in $x_h$ and $H^{s'}$ in $x_3$. Let $\langle x\rangle$ denote the quantity $\langle x\rangle\dd(1+|x|^2)^\frac12$, then for all $s,s' \in \mathbb R$, $H^{s, s'}$ is the space of tempered distributions $\psi \in \mathcal S'(\mathbb R^3)$ which satisfy
\[\nr \psi\nr_{{s,s'}}\dd\s  \langle\xi'\rangle^{2s}\,\langle \xi_3\rangle^{2s'}\, |\mathscr F \psi(\xi)|^2\,d\xi, \] 
where $\xi':=(\xi_1,\xi_2)$ and $\mathscr F$ denotes the Fourier transform. The space $\nr \psi\nr_{H^{s,s'}}$ endowed with the norm $\nr\cdot\nr_{s,s'}$ is a Hilbert space. Obviously, the homogenous anisotropic Sobolev space $\dot H^{s, s'}(\mathbb R^3)$ are obtained by replacing $\langle\cdot\rangle$ by $|\cdot|$. We will denote $L^p_h(L^q_v)$ the space $L^p(\mathbb R_{x_1}\times \mathbb R_{x_2};L^q(\mathbb R_{x_3}))$ endowed with the norm
\begin{align*}\nr \psi\nr_{L^p_h(L^q_v)}&\dd\nr\nr\psi\nr_{L^q(\mathbb R_{x_3})}\nr_{L^p(\mathbb R_{x_1}\times \mathbb R_{x_2})} = \left\{\int_{\mathbb R_{x_1}\times \mathbb R_{x_2}}  \left(\int_{\mathbb R_{x_3}} |\psi(x_h,x_3)|^q\,dx_3\right)^\frac pq\,dx_h \right\}^\frac1p.
\end{align*}
Equivalently, we denote $L^q_v(L^p_h)$ the space $L^q(\mathbb R_{x_3}(L^q({\mathbb R_{x_1}\times \mathbb R_{x_2}})))$ with the associated norm given by $\nr \psi\nr_{L^q_v(L^p_h)}:= \nr \nr \psi(\cdot,x_3)\nr_{L^p(\mathbb R_{x_1}\times \mathbb R_{x_2})}\nr_{L^q(\mathbb R_{x_3})}$.  The $L^p(\mathbb R^3)$ norms will be denoted $\nr\cdot\nr_p$. 
\vskip6pt
\noindent The mathematical analysis of the anisotropic Navier-Stokes system $\mathcal S_0$ was originally investigated in \cite{Chemin1} and \cite{Iftimie} where it is proved that the system $\mathcal S_0$ is locally well-posed for initial data in $H^{0,s}(\mathbb R^3)$ for all $s>\frac12$. Moreover, it has also been proved that if the initial data $u_0$ is such that 
\[\nr u_0\nr^{s-\frac12}_{L^2(\mathbb R^3)}\,\nr u_0\nr^{\frac32-s}_{\dot H^{0,s}(\mathbb R^3)} \leq c,\] 
for a sufficiently small constant $c$, then system $\mathcal S_0$ is globally well-posed. The aim of this short paper is to show how the damping term $|u|^{2\alpha}u$ gives rise to a smoothing effect in the vertical velocity.
Therefore it allows to get rid of the smallness assumption  (above) used in  $\mathcal S_0$. Even though, this result is still valid in the critical Sobolev and Besov spaces  $H^{0,\frac12}$ and $B^{0,\frac12}$ (see \cite{Paicu}), in order to avoid technicalities, in this paper we chose to  focus on a less optimal space to show how we take advantage of the damping term. A similar result in the spaces  $H^{0,\frac12}$ and $B^{0,\frac12}$  will be shown in a forthcoming paper soon.  Here, we specifically show the following  

\begin{thm}\label{thmnonoptimal}
Let $a,\nu>0, \alpha>1$ and $u_0 \in H^{0,1}(\mathbb R^3) $ such that ${\rm div} \,u_0=0$. Then,  system $\mathcal S_a$ has a unique global solution $u(t)$ satisfying 
\[u(t)\in L^\infty_{\rm loc}(\R;H^{0,1}(\mathbb R^3)) \cap L^2_{\rm loc}(\R; H^{1,1}(\RR))\cap L^{2\alpha+2}_{\rm loc}(\R; L^{2\alpha+2}(\RR)).\]
Moreover, the solution is in $C^0(\R; L^2(\mathbb R^3))$ and depends continuously on the initial data.
\end{thm}
\section*{Proof of Theorem \ref{thmnonoptimal}}
The rest of the paper is dedicated to the proof of Theorem \ref{thmnonoptimal}. In general, the proof is structured in four  steps. First, one defines a family of approximate systems $(\mathcal S_a^n)_{n\in \mathbb N}$ and show that this family has local in time smooth enough solutions $(u_n(t),p_n(t))$. This can be achieved for instance by the classical Friedrich's method. Second, one proves uniform bounds for  $(u_n(t),p_n(t))$ on some fixed time interval $[0,T]$. Next, one shows that the sequence of solutions to $(\mathcal S_a^n)_{n\in \mathbb N}$ converges towards some solution of $\mathcal S_a$ with adequate properties. Eventually, one exhibits a stability kind estimate leading to the continuous dependence of the solutions on the initial data, in particular their uniqueness. We refer to any textbook of fluid mechanics for technical details of this procedure (see e.g. \cite{Bahouri,Temam}). To shorten the presentation, we will only present the necessary uniform bounds by performing formal calculation using  system $\mathcal S_a$ instead of $(\mathcal S_a^n)_{n\in \mathbb N}$ and we will briefly explain how to pass to the limit.

\subsection*{ A priori Estimates}
We start by looking for an $L^2$ uniform estimate for the velocity. For this purpose, we multiply the first equation of system $\mathcal S_{a}$ by $u$ and integrate\footnote{This should be done on the smooth approximate solutions $u_n$.} over $\RR$ to get
\begin{align*}
\frac12\,\frac d{dt}\,\nr u(t)\nr_2^2 +\nu\,\nr \nabla_h u(t)\nr_2^2 + a\,\nr u(t)\nr_{2\alpha+2}^{2\alpha+2} =0,
\end{align*}
thanks to the fact that $\int_{\mathbb R^3} (u\cdot\nabla) u\cdot u\,dx=0$. Now, we integrate this equality with respect to time  
\begin{equation*}
\nr u(t)\nr_2^2 +2\nu\,\st\nr \nabla_h u(\tau)\nr_2^2\,d\tau + 2a\,\st\nr u(t)\nr_{2\alpha+2}^{2\alpha+2}\,d\tau= \nr u_0\nr_2^2.
\end{equation*}
This shows that if $u_0\in L^2(\mathbb R^3)$, then for all $t\in [0,T]$, it holds 
\begin{equation}\label{firstest}
u(t)\in L^\infty(\R;L^2(\mathbb R^3)) \cap L^2_{\rm loc}(\R;H^{1,0}(\mathbb R^3)) \cap L^{2\alpha+2}_{\rm loc}(\R;L^{2\alpha+2}(\mathbb R^3)).
\end{equation}
Next, we multiply the first equation of \sys\, by $-\partial_3^2 u$ and integrate\footnotemark[1]  \,over $\RR$  to get  
\begin{align}\label{brutres}
\frac12\,\frac d{dt}\,\nr \partial_3 u(t)\nr_2^2 +\nu\,\nr \nabla_h \partial_3 u(t)\nr_2^2 &\nonumber - \int_{\RR}\, (u(t)\cdot \nabla) u(t)\cdot \partial_3^2u(t)\,dx \\ &- a\int_\RR\, |u(t)|^{2\alpha}u(t)\cdot \partial_3^2u(t)\,dx =0.
\end{align}
Now, we handle the nonlinear terms. On the one hand, an integration by parts leads clearly to the fact that 
\begin{align}\label{good}
-\int_\RR\, |u|^{2\alpha}u\cdot \partial_3^2u\,dx&=\int_\RR\, |\partial_3 u|^2| u|^{2\alpha}\, dx +2\alpha  \int_\RR\, (u\cdot \partial_3)^2|u|^{2\alpha-2}\,dx \nonumber\\ &=(1+2\alpha)\,\nr |u|^\alpha\,\partial_3u\nr_2^2.
\end{align}
On the other hand, using integration by parts, we can write  
\begin{align*}
- \int_{\RR}\, (u\cdot \nabla) u\cdot \partial^2_3u\,dx &= \sum_{k,l=1}^3\,\int_\RR\partial_3 u_k\,\partial_k u_l\,\partial_3u_l\,dx \\
 &= \sum_{k=1}^2\sum_{l=1}^3\,\int_\RR\partial_3 u_k\,\partial_k u_l\,\partial_3u_l\,dx\\ &+\sum_{l=1}^3\,\int_\RR\partial_3 u_3\,\partial_3 u_l\,\partial_3u_l\,dx:=\mathscr T_1 +\mathscr T_2.
\end{align*}
Now, integrating again by parts, we obtain that
\begin{align*}
\mathscr T_1&= -\sum_{k=1}^2\sum_{l=1}^3 \int_\RR\,\left(u_l \,\partial_3u_l\, \partial_k\partial_3 u_k + u_l \,\partial_3u_k\, \partial_k\partial_3 u_l\right)\,dx.
\end{align*}
Moreover, using the fact that $\nabla\cdot u=0$, we have $-\partial_3u_3={\rm div}_h\,u_h$ where $u_h\dd (u_1,u_2)$. Thus
\begin{align*}
\mathscr T_2&= -\sum_{l=1}^3 \int_\RR\,{\rm div}_h\,u_h\,\partial_3 u_l\,\partial_3u_l \,dx.
\end{align*}
Next, using H\"older and Young inequality, it is rather easy to see that for all $f,g$ and $h$, we have for all $\alpha>1$ and $\epsilon_0,\epsilon_1>0$
\begin{align*}
\s\, f\,g\,h\,dx\leq \s\,|f|\,|g|^\frac1\alpha\,|g|^{1-\frac1\alpha}\,|h|\,dx 
&\leq \nr |f|\,|g|^\frac1\alpha\nr_{{2\alpha}}\,\nr |g|^{1-\frac1\alpha}|\nr_{\frac{2\alpha}{\alpha-1}}\,\nr h\nr_{2}\nonumber\\
&\leq \frac1{4\epsilon_0}\, \nr f^{\alpha}\,g\nr^{\frac2\alpha}_{2}\,\nr g\nr^{2(1-\frac1\alpha)}_{2}  +\epsilon_0\nr h\nr^2_{2}\nonumber\\
&\leq \frac{\epsilon_1}{4\epsilon_0\,}\, \nr f^{\alpha}\,g\nr^{2}_{2} + \frac{{\epsilon_1}^{\frac1{1-\alpha}}}{4\epsilon_0} \,\nr g\nr^{2}_{2}  +\epsilon_0\,\nr h\nr^2_{2}.
\end{align*}
Eventually, applying this inequality with $f=u_l, g=\partial_3 u_l, h= \partial_k\partial_3 u_k$  for the first part of $\mathscr T_1$, $f=u_l, g=\partial_3 u_k, h= \partial_k\partial_3 u_l$  for the second part of $\mathscr T_1$ and proceeding equivalently for $\mathscr T_2$, we obtain the existence of a constant $\gamma>0$ independent of $\alpha$ such that 
\begin{align}\label{mainineq}
 \int_{\RR}\, (u\cdot \nabla) u\cdot \partial^2_3u\,dx \leq \frac{\gamma\epsilon_1}{4\epsilon_0\,}\, \nr|u|^{\alpha}\,\partial_3\,u\nr^{2}_{2} + \frac{\gamma\epsilon_1^{\frac1{1-\alpha}}}{4\epsilon_0} \,\nr\partial_3\,u\nr^{2}_{2}  +\gamma\epsilon_0\,\nr \nabla_h\,\partial_3 u\nr^2_{2}.
 \end{align}
The idea then is to tune $\epsilon_0$ and $\epsilon_1$ to compensate the first and third terms of the right hand side of \eqref{mainineq} using \eqref{good} and the second term of the left hand side of \eqref{brutres}. More precisely, setting $\epsilon_0=\frac\nu{2\gamma}$ and $\epsilon_1=\frac{a\nu(1+4\alpha)}{\gamma^2}$ and using \eqref{good} and \eqref{mainineq}, the equality \eqref{brutres} implies the existence of some $\eta>0$ such that
\begin{align*}
\frac d{dt}\,\nr \partial_3 u(t)\nr_2^2 +\nu\,\nr \nabla_h \partial_3 u(t)\nr_2^2 +a\,\nr |u(t)|^\alpha\partial_3u(t)\nr_2^2  \leq\eta \,\nr \partial_3 u(t)\nr_2^2.
\end{align*}
Thus, thanks to Gronwall's inequality, we obtain the following bound
\[\nr \partial_3 u(t)\nr_2^2\leq \nr \partial_3 u_0\nr_2^2\,e^{\eta\,t},\quad \text{for all}\quad t\in [0,T].\]
 In particular, we have 
\begin{align*}
\nr \partial_3 u(t)\nr_2^2 +\nu\,\int_0^t\,\nr \nabla_h \partial_3 u(\tau)\nr_2^2\,d\tau +a\,\int_0^t\,\nr |u(\tau)|^\alpha\partial_3u(\tau)\nr_2^2\,d\tau  \leq \left(1+ e^{\eta\,t}\right)\nr \partial_3 u_0\nr_2^2.
\end{align*}
Thus, we infer
\begin{equation}\label{pass}
u(t)\in L^\infty(\R;H^{0,1}) \cap L^2(\R;H^{1,1}).
\end{equation}
Rigorously, these bounds hold for the approximate solutions constructed via the Friederich's regularization procedure. So, at this level, it remains only to pass to the limit in the sequence of solutions of  $(\mathcal S_a^n)_{n\in \mathbb N}$. For that purpose, the main point to show is that
\begin{align}\label{duality}\partial_t u&\in L^2_{\rm loc}(\R,H^{-1}(\RR)) + L^{1+\frac{1}{2\alpha+1}}_{\rm loc}(\R,L^{1+\frac{1}{2\alpha+1}}(\RR))\nonumber\\ &=\left(L^2_{\rm loc}(\R,H^1(\RR)) \cap L^{2\alpha+2}_{\rm loc}(\R,L^{2\alpha+2}(\RR))\right)^\star, \end{align}
where the star stands for the dual symbol. Indeed, let us recall the Lady\v{z}henskaya inequality, which is a special case of the Gagliardo-Nirenberg-Sobolev inequality (see e.g. \cite{Lady})
\begin{equation}\label{forduality1}
 \nr \psi\nr_4 \leq \delta_1 \,\nr \psi\nr_2^\frac14\,\nr \nabla\psi\nr_{2}^\frac34, \:\text{ for all }\: \psi\in  H^1_0(\RR).
\end{equation}
Therefore, using H\"older and Lady\v{z}henskaya inequalities, we have 
\begin{equation}\label{forduality2}\nr (u\cdot\nabla)u\nr_{{\frac43}} \leq \nr u\nr_4\,\nr \nabla u\nr_2 \leq \delta_1\nr u\nr_{2}^\frac14\,\nr  \nabla u\nr_{2}^\frac74 \leq \delta_1^8\,\nr u\nr_{2}^2 + \nr \nabla u \nr_{2}^2.
\end{equation}
Therefore $ (u\cdot\nabla)u \in L^2_{\rm loc}(\R,  H^{-1}(\RR))$. Also, we have clearly 
\begin{equation*}
\nr |u|^{2\alpha}u\nr^{1+\frac1{2\alpha+1}}_{1+\frac1{2\alpha+1}}=\nr u\nr^{2\alpha+2}_{2\alpha+2},\quad \text{and}\quad (2\alpha+2)^{-1} +(1+1/{2\alpha+1})^{-1}=1.
\end{equation*}
Thus, \eqref{duality} holds thanks to \eqref{firstest}, \eqref{pass}, \eqref{forduality1} and \eqref{forduality2}. Recall that 
\eqref{duality}  is needed in order to get some compactness in time. The passage to the limit follows using using classical argument by combining Ascoli's theorem and the Cantor diagonal process. 
\begin{remark}
It is rather standard to show that  $ u \in C^0(\R; L^2)$ by using \eqref{pass} and \eqref{duality}.  
 \end{remark}
\subsection*{ Uniqueness}
\noindent Now, we show the continuous dependence of the solutions on the initial data, in particular their uniqueness. Let $u(t)$ and $v(t)$ be two solutions of system \sys\, in the class $u(t)\in L^\infty(\R;H^{0,1}) \cap L^2(\R;H^{1,1})$. Let $w(t)=u(t)-v(t)$, then $w$ satisfies 
\begin{align*}
 \left\lbrace\begin{array}{rcll}
\partial_t\,w-\nu\,\Delta_h\,w +(w\cdot\nabla)\,u +(v\cdot\nabla)\,w+a\,|u|^{2\alpha}u-a\,|v|^{2\alpha}v &=&-\nabla\,(p_u-p_v),\\ \\
 \nabla\cdot w&=&0,\\ \\
u_{\vert_{t=0}}&=&u_0-v_0,& 
\end{array}
\right.
\end{align*}
We proceed as for the obtention of {\it a priori} estimates. Thanks to \eqref{duality}, the action of $\partial_t w$ on $w$ leads to 
\begin{align*}
\frac12\frac d{dt}\,\nr w\nr_2^2 + \nu\,\nr \nabla_hw\nr_2^2 +\int_{\mathbb R^3} (w\cdot \nabla) u \cdot w \,dx + a\,\int_{\mathbb R^3}\left(|u|^{2\alpha}u-\,|v|^{2\alpha}v\right)\,w\,dx =0.
\end{align*}
On the one hand 
\begin{align*}
\int_{\mathbb R^3} (w\cdot \nabla) u \cdot w \,dx =\sum_{k=1}^2 \sum_{l=1}^3 \,\int_{\mathbb R^3}\,w_k\,\partial_k u_l\,w_l\,dx+\sum_{l=1}^3 \,\int_{\mathbb R^3}\,w_3\,\partial_3 u_l\,w_l\,dx\dd\mathscr I_1 +\mathscr I_2.
\end{align*}
Now, using H\"older inequality, it holds
\begin{align*}
\mathscr I_1 &\leq \sum_{k=1}^2 \sum_{l=1}^3 \,\int _{\mathbb R} \,\nr \partial_k u_l\nr_{L^2_h}\,\nr w_k\nr_{L^4_h} \,\nr w_l\nr_{L^4_h} \,dx_3 \\ &\leq  \sum_{k=1}^2 \sum_{l=1}^3 \,\nr \partial_k u_l\nr_{L^\infty_v(L^2_h)}\,\nr w_k\nr_{L^2_v(L^4_h)} \,\nr w_l\nr_{L^2_v(L^4_h)}.
\end{align*}
Now, using the Sobolev embedding $\dot H^\frac12_h \hookrightarrow L^4_h$ and interpolating $\dot H^\frac12_h$ between $\dot H^1_h$ and $L^2_h$, we obtain clearly  for all $\psi \in L^2_v \cap \dot H^1_h$
\begin{align}\label{p1}
\nr \psi\nr_{L^2_v(L^4_h)} \leq C\,\nr \nabla_h\psi\nr^\frac12_{2}\,\nr \psi\nr^\frac12_{2}.
\end{align}
Also, we have 
\begin{align}\label{p2}
\nr \psi(\cdot, x_{3})\nr^2_{L^2_h} &= \int_{-\infty}^{x_3} \,\frac d{dz} \left(\nr \psi(\cdot,z)\nr^2_{L^2_h}\right) \,dz\nonumber\\
&= 2 \int_{-\infty}^{x_3} \int_{\mathbb R^2}\, \psi(x_h,z)\,\partial_z\psi(x_h,z) dx_h\,dz \leq 2\,\nr \psi\nr_2\,\nr\partial_3\psi\nr_2.
\end{align}
Therefore, using \eqref{p1} with $\psi=w_k$ and $\psi=w_l$ and \eqref{p2} with $\psi=\partial_k u_l$, we obtain thanks to Young's inequality
\begin{align*} 
\mathscr I_1&\leq C\,\nr \partial_3\nabla_h u\nr^\frac12_{2}\,\nr \nabla_h u\nr^\frac12_{2}\,\nr \nabla_h w\nr_{2}\,\nr w\nr_{2} 
\leq \frac\nu4\,\nr \nabla_h w\nr_{2}^2 + {C}\, \left(\nr \partial_3\nabla_h u\nr^2_{2}+\nr \nabla_h u\nr^2_{2} \right)\,\nr w\nr^2_{2} .
\end{align*}
Next, proceeding in the same way, we get  
\begin{align*} 
\mathscr I_2&\leq C\, \nr w_3\nr_{L^\infty_v(L^2_h)}\,\nr \partial_3\nabla_h u\nr^\frac12_{2}\,\nr \partial_3 u\nr^\frac12_{2}\,\nr \nabla_h w\nr^\frac12_{2}\,\nr w\nr^\frac12_{2}.
\end{align*}
But, using the fact that $\nabla\cdot w=0$, thus ${\rm div}_h w_h =-\partial_3 w_3$ , we get 
\begin{align*}
\nr w_3\nr^2_{L^2_h} &=  2 \int_{-\infty}^{x_3} \int_{\mathbb R^2}\, w_3(x_h,z)\,\partial_3w_3(x_h,z) dx_h\,dz \\ &=  -2 \int_{-\infty}^{x_3} \int_{\mathbb R^2}\, w_3(x_h,z)\,{\rm div}_h w_h(x_h,z) dx_h\,dz \\
&\leq 2\,\nr {\rm div}_h w_h\nr_2\,\nr w_3\nr_2.
\end{align*}
Hence
\begin{align*}
\mathscr I_2&\leq C\, \nr {\rm div}_h w_h\nr^\frac12_2\,\nr w_3\nr^\frac12_2\,\nr \partial_3\nabla_h u\nr^\frac12_{2}\,\nr \partial_3 u\nr^\frac12_{2}\,\nr \nabla_h w\nr^\frac12_{2}\,\nr w\nr^\frac12_{2}\\
&\leq \frac\nu4\,\nr \nabla_h w\nr_{2}^2 + C\,\left(\nr \partial_3\nabla_h u\nr^2_{2} +\nr\partial_3 u\nr^2_{2} \right)\,\nr w\nr_2^2
\end{align*}
On the other hand, it is well known that there exists a nonnegative constant ${\kappa}={\kappa}(\alpha)$ such that
\begin{equation*}
0\leq{\kappa}\,|u-v|^{2}\,\left(|u|+|v|\right)^{2\alpha}\leq \left(|u|^{2\alpha}u-|v|^{2\alpha}v\right)\cdot (u-v).
\end{equation*} 
Thus
\begin{align*}
a\,\int_{\mathbb R^3}\left(|u|^{2\alpha}u-\,|v|^{2\alpha}v\right)\,w\,dx \geq a\,{\kappa}\,\int_{\mathbb R^3}\left(|u|+|v|\right)^{2\alpha}\,w^2\,dx = a{\kappa}\,\nr \left(|u|+|v|\right)^{\alpha}\,w\nr^2_2.
\end{align*}
All in all, we have 
\begin{align*}
\frac d{dt}\,\nr w\nr_2^2 + \nu\,\nr \nabla_hw\nr_2^2 &+2a{\kappa}\,\nr \left(|u|+|v|\right)^{\alpha}\,w\nr^2_2 \\& \leq  C\,\left(\nr \partial_3\nabla_h u\nr^2_{2} +\nr\partial_3 u\nr^2_{2} + \nr \nabla_h u\nr^2_{2}\right)\,\nr w\nr_2^2
\end{align*}
Setting  $\mathcal L(t):= C\,\left(\nr \partial_3\nabla_h u\nr^2_{2} +\nr\partial_3 u\nr^2_{2} +\nr \nabla_h u\nr^2_{2}\right)$ and integrating the above inequality with respect to time, we obtain that for all $t\in [0,T]$
\begin{align}\label{togron}
\nr w(t)\nr_2^2 +\nu\, \int_0^t\,\nr \nabla_hw(\tau)\nr_2^2\,d\tau &+2a{\kappa}\, \int_0^t\,\nr \left(|u|+|v|\right)^{\alpha}\,w(\tau)\nr^2_2\,d\tau \nonumber\\  &\leq \nr w_0\nr_2^2 + \int_0^t\,\mathcal L(\tau)\,\nr w(\tau)\nr_2^2\,d\tau
\end{align}
Now, since $\partial_3 u\in L^\infty_{\rm loc}(\R,L^2(\RR))$, $\partial_3\nabla_h u\in L^2_{\rm loc}(\R, L^2(\RR))$
and $\nabla_h u\in L^2_{\rm loc}(\R, L^2(\RR))$, it is clear that $\mathcal L(t)\in L^1_{\rm loc}(\R)$. Therefore, Gronwall's Lemma applied to the inequality \eqref{togron} leads to the uniqueness for all $\alpha>1$ and the proof is complete.
\vspace{0.5cm}

\textbf{Acknowledgment:} The research of Hakima Bessaih was partially supported by NSF grant DMS-1418838. The research of  Saber Trabelsi was supported by the King Abdullah University of Science and Technology.


\end{document}